\documentclass{amsart}

\theoremstyle{definition}

\theoremstyle{remark}

\usepackage{graphicx}
\usepackage{color}
\usepackage{amsmath}
\usepackage{amsfonts}
\usepackage{amssymb}
\newcommand{\be}{\begin{equation}}
\newcommand{\ee}{\end{equation}}

\newcommand{\ba}{\begin{array}}

\newcommand{\ea}{\end{array}}

\newcommand{\beq}{\begin{eqnarray}}

\newcommand{\eeq}{\end{eqnarray}}

\newtheorem{lm}{lemma}

\newtheorem{thee}{theorem}

\newtheorem{proo}{proposition}

\newtheorem{co}{corollary}

\newtheorem{rem}{remark}

\newtheorem{deff}{definition}

\newcommand{\bd}{\begin{deff}}

\newcommand{\ed}{\end{deff}}

\newcommand{\bl}{\begin{lm}}

\newcommand{\el}{\end{lm}}

\newcommand{\bp}{\begin{proo}}

\newcommand{\ep}{\end{proo}}

\newcommand{\bt}{\begin{thee}}

\newcommand{\et}{\end{thee}}

\newcommand{\bc}{\begin{co}}

\newcommand{\ec}{\end{co}}

\newcommand{\brm}{\begin{rem}}

\newcommand{\erm}{\end{rem}}

\newcommand{\der}{{\rm d}}

\usepackage{t1enc}

\newcommand{\newc}{\newcommand}

\let\ccdot\cdot
\def\cdot{\hbox to 2.5pt{\hss$\ccdot$\hss}}

\newc{\aR}{\mbox{\boldmath{$ R$}}}
\newc{\aS}{\mbox{\boldmath{$ S$}}}
\newc{\aT}{\mbox{\boldmath{$ T$}}}
\newc{\aW}{\mbox{\boldmath{$ W$}}}

\newc{\aK}{\mbox{\boldmath{$ K$}}}
\newc{\aL}{\mbox{\boldmath{$ L$}}}

\newcommand{\bbC}{\mathbb{C}}

\usepackage{amssymb}
\usepackage{amscd}

\newc{\obstrn}[2]{B^{#1}_{#2}}

\newcommand{\rpl}                         
{\mbox{$
\begin{picture}(12.7,8)(-.5,-1)
\put(0,0.2){$+$}
\put(4.2,2.8){\oval(8,8)[r]}
\end{picture}$}}

\newcommand{\lpl}                         
{\mbox{$
\begin{picture}(12.7,8)(-.5,-1)
\put(2,0.2){$+$}
\put(6.2,2.8){\oval(8,8)[l]}
\end{picture}$}}

\usepackage{ifthen}

\newcommand{\bbR}{\mathbb{R}}

\newc{\tensor}[1]{#1}
\newc{\Mvariable}[1]{\mbox{#1}}
\newc{\down}[1]{{}_{#1}}
\newc{\up}[1]{{}^{#1}}

\newc{\JulyStrut}{\rule{0mm}{6mm}}
\newc{\midtenPan}{\mbox{\sf S}}
\newc{\midten}{\mbox{\sf T}}
\newc{\midtenEi}{\mbox{\sf U}}
\newc{\ATen}{\mbox{\sf E}}
\newc{\BTen}{\mbox{\sf F}}
\newc{\CTen}{\mbox{\sf G}}

\def\sideremark#1{\ifvmode\leavevmode\fi\vadjust{\vbox to0pt{\vss
 \hbox to 0pt{\hskip\hsize\hskip1em
 \vbox{\hsize3cm\tiny\raggedright\pretolerance10000
 \noindent #1\hfill}\hss}\vbox to8pt{\vfil}\vss}}}%

\numberwithin{equation}{section}

\newcounter{romenumi}
\newcommand{\labelromenumi}{(\roman{romenumi})}

\begin{document}
\title{Construction of conjugate functions}

\author{Pawe\l~ Nurowski} \address{Instytut Fizyki Teoretycznej,
Uniwersytet Warszawski, ul. Hoza 69, Warszawa, Poland}
\email{nurowski@fuw.edu.pl} \thanks{This research was supported by
the KBN grant 1 P03B 07529}

\date{\today}

\begin{abstract} We find all pairs of real analytic 
functions $f$ and $g$ in $\bbR^n$ such that $|\nabla f|=|\nabla g|$ and 
$(\nabla f)(\nabla g)=0$.

\end{abstract}
\maketitle
We consider the following problem:
Let $n\geq 3$.
Find all pairs of functions $(f,g)$, $f:\bbR^n\to \bbR$, $g: \bbR^n\to \bbR$ such
that
$$|\nabla f|=|\nabla g|,\quad\quad\quad (\nabla f)(\nabla g)=0.$$ 
Functions $f$ and $g$ constituting such a pair 
are called {\it conjugate}.

Finding a pair of conjugate functions $(f,g)$ is equivalent to finding a 
complex valued function $h:\bbR^n\to \bbC$ such that
$$(\nabla h)(\nabla h)=0.$$
Having $h$ one gets $f$ and $g$ as its real and imaginary parts,
respectively.

Solving the problem at algebraic level, one needs to find a nice representation
for {\it null} vectors in $\bbC^n$. This may be done by means of spinors.

The solution of the problem in dimension $n=3$ is as follows.

We consider $\bbR^3$ with coordinates $(x,y,z)$. The most general form of $\nabla h$ is given in terms of two
complex-valued functions (a spinor) $(\phi_1,\phi_2)$ 
$$h_x=\phi_1^2-\phi_2^2,\quad\quad\quad
h_y=i(\phi_1^2+\phi_2^2),\quad\quad\quad h_z=2\phi_1\phi_2.$$ 

Now, the integrability condition for existence of $h$ is 
$$d~[~(\phi_1^2-\phi_2^2)d x+i(\phi_1^2+\phi_2^2)d y+2\phi_1\phi_2 dz~]=0,$$
\bigskip
which is equivalent to 
$$[\phi_1(dx+idy)+\phi_2dz]\wedge d\phi_1+[-\phi_2(dx-idy)+\phi_1
dz]\wedge d\phi_2=0.$$
This motivates introduction of two functions
$$ X_1=\phi_1 (x+iy)+\phi_2 z,\quad\quad\quad X_2=-\phi_2(x-iy)+\phi_1 z.$$
Having them the integrability condition is:
$$\der(X_1\der\phi_1+X_2\der\phi_2)=0.$$
Its general solution analytic in $(x,y,z)$ is
$$X_1=F_1(\phi_1,\phi_2)\quad\quad\quad X_2=F_2(\phi_1,\phi_2),$$
where $F=F(\phi_1,\phi_2)$ is a complex valued function analytic in
both variables, and $F_1=\tfrac{\partial F}{\partial \phi_1}$, $F_2=\tfrac{\partial F}{\partial \phi_2}$.

Explicitely, one finds $\phi_1$ and $\phi_2$ specifying $F$ and
solving the algebraic equations
$$\phi_1(x+iy)+\phi_2z=F_1(\phi_1,\phi_2)$$
$$-\phi_2(x-iy)+\phi_1 z=F_2(\phi_1,\phi_2).$$

Once $(\phi_1,\phi_2)$ is found then $h=\int (\phi_1^2-\phi_2^2)dx$,
which solves the problem. 

So in dimensin $n=3$ the conclusion is that the general analytic solution is
generetated by one complex analytic function $F$ of two variables.

Although, for arbitrary $n$ we can not use spinors anymore, a similar 
procedure works for any $n\geq 3$.

In $\bbR^n$ we introduce coordinates $(x_1,x_2,...,x_n)$ and parametrize 
a general null vector $(h_{x_1},h_{x_2},...,h_{x_n})$ by means of $n-1$ 
complex functions $(\phi_1,\phi_2,...,\phi_{n-1})$ via 
\begin{eqnarray*}
h_{x_1}&=&\phi_1^2+\phi_2^2+...+\phi_{n-2}^2-\phi_{n-1}^2\\
h_{x_2}&=&i(\phi_1^2+\phi_2^2+...+\phi_{n-2}^2+\phi_{n-1}^2)\\
h_{x_3}&=&2\phi_1\phi_{n-1}\\
...&&\\
h_{x_k}&=&2\phi_{k-2}\phi_{n-1}\\
...&&\\
h_{x_n}&=&2\phi_{n-2}\phi_{n-1}.
\end{eqnarray*}
This parametrization is known in the theory od minimal surfaces \cite{wood}.

Now the integrability condition
$$\der~[h_{x_1}\der x_1+h_{x_2}\der x_+...+h_{x_n}\der x_n]=0$$
for existence of $h:\bbR^n\to \bbC$ is
$$\der[X_1\der\phi_1+X_2\der\phi_2+...+X_{n-1}\der\phi_{n-1}]=0$$
with $n-1$ functions $(X_1,X_2,...,X_{n-1})$ given by
\begin{eqnarray*}
X_1&=&\phi_1(x_1+ix_2)+\phi_{n-1}x_3\\
X_2&=&\phi_2(x_1+ix_2)+\phi_{n-1}x_4\\
...&&\\
X_{n-2}&=&\phi_{n-2}(x_1+ix_2)+\phi_{n-1}x_n\\
X_{n-1}&=&-\phi_{n-1}(x_1-ix_2)+\phi_1x_3+\phi_2x_4+...+\phi_{n-2}x_n.
\end{eqnarray*}
Thus, similarly to the $n=3$ case, we first choose an arbitrary holomorphic function $F=F(\phi_1,\phi_2,...,\phi_{n-1})$ of $n-1$ complex variables and, denoting its derivatives by $F_i=\tfrac{\partial F}{\partial\phi_i}$, $i=1,2,...,n-1$, solve algebraic equations 
\begin{eqnarray*}
F_1&=&\phi_1(x_1+ix_2)+\phi_{n-1}x_3\nonumber\\
F_2&=&\phi_2(x_1+ix_2)+\phi_{n-1}x_4\nonumber\\
...&&\label{eq1}\\
F_{n-2}&=&\phi_{n-2}(x_1+ix_2)+\phi_{n-1}x_n\nonumber\\
F_{n-1}&=&-\phi_{n-1}(x_1-ix_2)+\phi_1x_3+\phi_2x_4+...+\phi_{n-2}x_n,\nonumber
\end{eqnarray*}
for $(\phi_1,\phi_2,...,\phi_{n-1})$ as functions of $(x_1,x_2,...,x_n)$.
Then $h$ is found by simple integration to be e.g. 
$h=\int 2\phi_{n-2}\phi_{n-1}\der x_n$.

We note that, under some regularity assumptions on the function $F$, 
equations (\ref{eq1}) associate to any solution $h$ an 
$n$-dimensional real hypersurface $M_n$ embedded in $\bbC^{n-1}$. One gets the equations for this hypersurface in coordinates 
$(\phi_1,\phi_2,...,\phi_{n-1})\in\bbC^{n-1}$ by elliminating the real parameters  
$(x_1,x_2,...,x_n)$ from equations (\ref{eq1}). For example if $n=3$, given $F=F(\phi_1,\phi_2)$, we have a 3-dimensional real hypersurface $M_3$ in $\bbC^2$ defined by
$$M_3=\{~(\phi_1,\phi_2)\in\bbC^2~:~ {\rm Im}(F_1\bar{\phi}_2-F_2\bar{\phi}_1)=0~\}.$$
In case of $n=4$, given $F=F(\phi_1,\phi_2,\phi_3)$ we have
\begin{eqnarray*}
M_4&=&\{~(\phi_1,\phi_2,\phi_3)\in\bbC^3~:~\\
0&=&{\rm Im}\big(~\bar{\phi}_3~[F_1(\phi_2^2\bar{\phi}_1^2+|\phi_2|^4-|\phi_1\phi_3|^2-|\phi_3|^4)-\\
&&F_2(|\phi_2|^2\phi_1\bar{\phi}_2+(|\phi_1|^2+|\phi_3|^2)\phi_2\bar{\phi}_1)+F_3\phi_3(\phi_1\bar{\phi}_2^2+(|\phi_1|^2+|\phi_3|^2)\bar{\phi}_1)]~\big)\\
0&=&{\rm Im}\big(~\bar{\phi}_3~[F_1(|\phi_1|^2\phi_2\bar{\phi}_1+(|\phi_2|^2+|\phi_3|^2)\phi_1\bar{\phi}_2)-\\
&&F_2(\phi_1^2\bar{\phi}_2^2+|\phi_1|^4-|\phi_2\phi_3|^2-|\phi_3|^4)-F_3\phi_3(\phi_2\bar{\phi}_1^2+(|\phi_2|^2+|\phi_3|^2)\bar{\phi}_2)]~\big)~\}.
\end{eqnarray*}
We also note that hypersurfaces $M_n$ are foliated by $(n-2)$-dimensional leaves which are the images under the map $(x_1,x_2,...,x_n)\mapsto(\phi_1,\phi_2,...,\phi_{n-1})$ of the intersections in $\bbR^n$ of the level surfaces  $f=c_1$ and $g=c_2$ corresponding to the conjugate functions $(f,g)$ associated with $h$. Thus, in particular, $M_3$ has a distinguished foliation by real curves and $M_4$ has a distinguished foliation by real surfaces.

Returning to our original problem it is interesting to note 
that slightly different procedure, based on a {\it three}-linear representation 
of a complex null vector, can be used to solve the problem in dimension $n=5$.

In $\bbR^5$ we use coordinates $(x,y,z,t,u)$.
Now we write a null vector $(h_x,h_y,h_z,h_t,h_u)$ in terms of {\it six} functions $(\phi_1,\phi_2,\phi_3,\phi_4,\phi_5,\phi_6)$ as follows:
$$h_x=i(\phi_1\phi_2 + \phi_3\phi_4)\phi_5 - 
    \tfrac12(\phi_1^2 + \phi_2^2 - \phi_3^2 - \phi_4^2)\phi_6,$$
$$
h_y=i (-\phi_1 \phi_3  + \phi_2 \phi_4)\phi_5 + (\phi_2 \phi_3 - \phi_1 \phi_4)\phi_6$$
$$
h_z=(-\phi_1\phi_2+ \phi_3\phi_4)\phi_5 -\tfrac{i}{2}(\phi_1^2 + \phi_2^2 + \phi_3^2 +\phi_4^2)\phi_6
$$
$$
h_t=(\phi_1 \phi_3 + \phi_2 \phi_4) \phi_5
$$
$$
h_u=(\phi_1 \phi_3 + \phi_2 \phi_4) \phi_6.
$$
Check that $h_x^2+h_y^2+h_z^2+h_t^2+h_u^2=0$. 
Remarkably again, the integrability conditions
$$\der~[h_x\der x+h_y\der y+h_z\der z+h_t \der t+h_u\der u~]~=0$$
are equivalent to
$$\der(X\der \phi_1+Y\der\phi_2+Z\der\phi_3+T\der\phi_4+U\der\phi_5+W\der \phi_6)=0,
$$
where
$$X=\phi_3\phi_5 (t-i y)  -( \phi_1 \phi_6- i \phi_2 \phi_5) (x+iz)  - \phi_4 \phi_6 y + \phi_3\phi_6 u $$
$$Y=\phi_4 \phi_5 (t +iy) - (\phi_2 \phi_6 -i \phi_1 \phi_5)(x+iz) + 
  \phi_3 \phi_6 y+ \phi_4 \phi_6 u$$
$$Z=\phi_1 \phi_5 (t-iy)   + (\phi_3 \phi_6 + i \phi_4 \phi_5) (x-iz) + \phi_2 \phi_6 y + \phi_1 \phi_6 u$$
$$T=\phi_2 \phi_5( t+iy) + (\phi_4 \phi_6+ i \phi_3 \phi_5) (x-iz) - 
  \phi_1 \phi_6 y  + \phi_2 \phi_6 u $$
$$U=\phi_1 \phi_3 (t-iy) + 
  \phi_2 \phi_4 (t+iy) + i \phi_1 \phi_2 (x+iz) + i \phi_3 \phi_4 (x-iz) $$
$$W=-\tfrac12 (\phi_1^2 + \phi_2^2) (x +iz)+\tfrac12  
        (\phi_3^2 + \phi_4^2 )(x  - i z)+ ( \phi_2 \phi_3   - 
        \phi_1\phi_4) y+(\phi_1\phi_3  + \phi_2\phi_4) u.$$
Thus taking an analytic function $F=F(\phi_1,\phi_2,\phi_3,\phi_4,\phi_5,\phi_6)$ of six variables, and solving algebraic equations
$$ F_1=\phi_3\phi_5 (t-i y)  -( \phi_1 \phi_6- i \phi_2 \phi_5) (x+iz)  - \phi_4 \phi_6 y + \phi_3\phi_6 u$$
$$F_2=\phi_4 \phi_5 (t +iy) - (\phi_2 \phi_6 -i \phi_1 \phi_5)(x+iz) + 
  \phi_3 \phi_6 y+ \phi_4 \phi_6 u)$$
$$F_3=\phi_1 \phi_5 (t-iy)   + (\phi_3 \phi_6 + i \phi_4 \phi_5) (x-iz) + \phi_2 \phi_6 y + \phi_1 \phi_6 u$$
$$F_4=\phi_2 \phi_5( t+iy) + (\phi_4 \phi_6+ i \phi_3 \phi_5) (x-iz) - 
  \phi_1 \phi_6 y  + \phi_2 \phi_6 u $$
$$F_5=\phi_1 \phi_3 (t-iy) + 
  \phi_2 \phi_4 (t+iy) + i \phi_1 \phi_2 (x+iz) + i \phi_3 \phi_4 (x-iz) $$
$$F_6=-\tfrac12 (\phi_1^2 + \phi_2^2) (x +iz)+\tfrac12  
        (\phi_3^2 + \phi_4^2 )(x  - i z)+ ( \phi_2 \phi_3   - 
        \phi_1\phi_4) y+(\phi_1\phi_3  + \phi_2\phi_4) u$$
for $(\phi_1,\phi_2,\phi_3,\phi_4,\phi_5,\phi_6)$ as functions of $(x,y,z,t,u)$, we generate a solution to $(\nabla h)(\nabla h)=0$ via, for example, an integral $h=\int[(\phi_1 \phi_3 + \phi_2 \phi_4) \phi_6]\der u$.

This note is motivated by a paper \cite{east}. It reminds very much the Kerr theorem (see e.g. \cite{pn,at}).

\end{document}